\documentclass[11pt]{amsart}

\begin{document}

\def\dbl{[\hskip -1pt[}
\def\dbr{]\hskip -1pt]}

\title{Analytic regularity of CR maps into spheres}
\author{Nordine Mir}
\address{Universit\'e de Rouen, Laboratoire de Math\'ematiques Rapha\"el Salem, UMR 6085 CNRS, 76821 Mont-Saint-Aignan Cedex, France}
\email{Nordine.Mir@univ-rouen.fr}
\thanks{\noindent 2000 {{\em Mathematics Subject Classification.}  32H02, 32H04, 32V20, 32V30,
32V40\\
Math. Res. Lett. {\bf 10}, no.\ 4, (2003),
447--457.}}

\def\Label#1{\label{#1}}
\def\1#1{\ov{#1}}
\def\2#1{\widetilde{#1}}
\def\6#1{\mathcal{#1}}
\def\4#1{\mathbb{#1}}
\def\3#1{\widehat{#1}}
\def\K{{\4K}}
\def\LL{{\4L}}

\def\C{{\4C}}
\def\R{{\4R}}
\def \MM{{\4M}}
\def \S{{\4S}^{2N'-1}}

\def \B{{\4B}^{2N'-1}}

\def \H{{\4H}^{2l-1}}

\def \F{{\4H}^{2N'-1}}

\def \LL{{\4L}}

\def\Re{{\sf Re}\,}
\def\Im{{\sf Im}\,}

\def\s{s}
\def\k{\kappa}
\def\ov{\overline}
\def\span{\text{\rm span}}
\def\ad{\text{\rm ad }}
\def\tr{\text{\rm tr}}
\def\xo {{x_0}}
\def\Rk{\text{\rm Rk\,}}
\def\sg{\sigma}
\def \emxy{E_{(M,M')}(X,Y)}
\def \semxy{\scrE_{(M,M')}(X,Y)}
\def \jkxy {J^k(X,Y)}
\def \gkxy {G^k(X,Y)}
\def \exy {E(X,Y)}
\def \sexy{\scrE(X,Y)}
\def \hn {holomorphically nondegenerate}
\def\hyp{hypersurface}
\def\prt#1{{\partial \over\partial #1}}
\def\det{{\text{\rm det}}}
\def\wob{{w\over B(z)}}
\def\co{\chi_1}
\def\po{p_0}
\def\fb {\bar f}
\def\gb {\bar g}
\def\Fb {\ov F}
\def\Gb {\ov G}
\def\Hb {\ov H}
\def\zb {\bar z}
\def\wb {\bar w}
\def \qb {\bar Q}
\def \t {\tau}
\def\z{\chi}
\def\w{\tau}
\def\Z{\zeta}
\def\phi{\varphi}
\def\eps{\varepsilon}

\def \T {\theta}
\def \Th {\Theta}
\def \L {\Lambda}
\def\b {\beta}
\def\a {\alpha}
\def\o {\omega}
\def\l {\lambda}

\def \im{\text{\rm Im }}
\def \re{\text{\rm Re }}
\def \Char{\text{\rm Char }}
\def \supp{\text{\rm supp }}
\def \codim{\text{\rm codim }}
\def \Ht{\text{\rm ht }}
\def \Dt{\text{\rm dt }}
\def \hO{\widehat{\mathcal O}}
\def \cl{\text{\rm cl }}
\def \bR{\mathbb R}
\def \bS{\mathbb S}
\def \bK{\mathbb K}
\def \bD{\mathbb D}
\def \bC{\mathbb C}
\def \C{\mathbb C}
\def \N{\mathbb N}
\def \bL{\mathbb L}
\def \bZ{\mathbb Z}
\def \bN{\mathbb N}
\def \scrF{\mathcal F}
\def \scrK{\mathcal K}
\def \mc #1 {\mathcal {#1}}
\def \scrM{\mathcal M}
\def \cR{\mathcal R}
\def \scrJ{\mathcal J}
\def \scrA{\mathcal A}
\def \scrO{\mathcal O}
\def \scrV{\mathcal V}
\def \scrL{\mathcal L}
\def \scrE{\mathcal E}
\def \hol{\text{\rm hol}}
\def \aut{\text{\rm aut}}
\def \Aut{\text{\rm Aut}}
\def \J{\text{\rm Jac}}
\def\jet#1#2{J^{#1}_{#2}}
\def\gp#1{G^{#1}}
\def\gpo{\gp {2k_0}_0}
\def\emmp {\scrF(M,p;M',p')}
\def\rk{\text{\rm rk\,}}
\def\Orb{\text{\rm Orb\,}}
\def\Exp{\text{\rm Exp\,}}
\def\Span{\text{\rm span\,}}
\def\d{\partial}
\def\D{\3J}
\def\pr{{\rm pr}}

\def \CZZ {\C \dbl Z,\zeta \dbr}
\def \D{\text{\rm Der}\,}
\def \Rk{\text{\rm Rk}\,}
\def \CR{\text{\rm CR}}
\def \ima{\text{\rm im}\,}
\def \I {\mathcal I}

\def \M {\mathcal M}

\newtheorem{Thm}{Theorem}[section]
\newtheorem{Cor}[Thm]{Corollary}
\newtheorem{Pro}[Thm]{Proposition}
\newtheorem{Lem}[Thm]{Lemma}

\theoremstyle{definition}\newtheorem{Def}[Thm]{Definition}

\theoremstyle{remark}
\newtheorem{Rem}[Thm]{Remark}
\newtheorem{Exa}[Thm]{Example}
\newtheorem{Exs}[Thm]{Examples}

\numberwithin{equation}{subsection}

\def\bl{\begin{Lem}}
\def\el{\end{Lem}}
\def\bp{\begin{Pro}}
\def\ep{\end{Pro}}
\def\bt{\begin{Thm}}
\def\et{\end{Thm}}
\def\bc{\begin{Cor}}
\def\ec{\end{Cor}}
\def\bd{\begin{Def}}
\def\ed{\end{Def}}
\def\br{\begin{Rem}}
\def\er{\end{Rem}}
\def\be{\begin{Exa}}
\def\ee{\end{Exa}}
\def\bpf{\begin{proof}}
\def\epf{\end{proof}}
\def\ben{\begin{enumerate}}
\def\een{\end{enumerate}}

\keywords{}

\maketitle

\begin{abstract}
Let $M\subset \C^N$ be a connected real-analytic hypersurface and
$\S\subset \C^{N'}$ the unit real sphere, $N'> N\geq 2$. Assume
that $M$ does not contain any complex-analytic hypersurface of
$\C^N$ and that there exists at least one strongly pseudoconvex
point on $M$. We show that any CR map $f\colon M\to \S$ of class
$\6C^{N'-N+1}$ extends holomorphically to a neighborhood of $M$ in
$\C^N$.
\end{abstract}


\section{Introduction}\Label{int}

In this paper we are interested in the analytic regularity of CR
mappings from real-analytic hypersurfaces into higher dimensional
unit spheres in complex space. While there is a wide literature
deciding when CR maps, of a given smoothness, between two
real-analytic hypersurfaces in the same complex space must be
real-analytic (see e.g.\ \cite{BN, Fo2, BERbook, HH} for complete
references up to 1999), very little is known about the analyticity
of such maps when the hypersurfaces lie in complex spaces of
different dimension. The case of CR maps with target unit spheres,
arising e.g.\ from the embedding problem for pseudoconvex domains
into balls (see e.g.\ \cite{F86, EHZ}), has attracted the
attention of many authors. The first regularity result in such a
situation was given by Webster \cite{W79} who showed that any CR
map of class $\6C^3$ from a real-analytic strongly pseudoconvex
hypersurface in $\C^N$ into the unit sphere ${{\4S}^{2N+1}}\subset
\C^{N+1}$ is real-analytic on a dense open subset of the source
hypersurface. Later Forstneri\v c \cite{F} generalized Webster's
result by showing that the same conclusion holds for any CR map of
class $\6C^{N'-N+1}$ with an arbitrary unit sphere $\S\subset
\C^{N'}$ as a target, $N'> N\geq 2$ (see also \cite{Hu94}). He
also asked in the same paper whether the real-analyticity, or
equivalently, the holomorphic extendability of such maps holds at
every point. A partial positive answer in codimension one (i.e.\
for $N'-N=1$) was given by Baouendi, Huang and Rothschild
\cite{BHR} for the case of CR maps from a real-algebraic
hypersurface of D'Angelo finite type (i.e.\ not containing any
positive dimensional complex-analytic subvariety) in $\C^N$ into
the unit sphere ${{\4S}^{2N+1}}\subset \C^{N+1}$. (In \cite{BHR},
the required smoothness for the maps depends on the so-called
D'Angelo type of the reference point, which is always greater or
equal to two.) In this paper we prove the following theorem, which
provides a positive answer to Forstneri\v c's question and at the
same time generalizes the quite independent result of
Baouendi-Huang-Rothschild.

\begin{Thm}\Label{main1}
Let $M\subset \C^N$ be a connected real-analytic hypersurface and
$\S\subset \C^{N'}$ the unit sphere, $N'> N\geq 2$. Assume that
$M$ does not contain any complex-analytic hypersurface of $\C^N$
and that there exists at least one strongly pseudoconvex point on
$M$. Then any CR map $f\colon M\to \S$ of class $\6C^{N'-N+1}$
extends holomorphically to a neighborhood of $M$ in $\C^N$.
\end{Thm}

Observe that if $M$ is a real hypersurface as in Theorem
\ref{main1} and admits a non-constant CR map into a higher
dimensional unit sphere, then it has to be pseudoconvex and
therefore the set of strongly pseudoconvex points forms a Zariski
open subset of $M$ (see Lemmas \ref{hor1} and \ref{hor2}). Hence
Theorem \ref{main1} gives the holomorphic extension of $f$ not
only near any strongly pseudoconvex point of $M$, but also near
all weakly pseudoconvex ones. It is worth noticing that, under the
assumptions of Theorem \ref{main1}, the set of weakly pseudoconvex
points of $M$ may even contain complex curves. Example \ref{add}
below provides a simple illustration of such a situation where
Theorem \ref{main1} applies. On the other hand, Theorem \ref{main1}
in conjunction with Lemma \ref{hor1} (ii) below implies the following
corollary for real-analytic
hypersurfaces of D'Angelo finite type.

\begin{Cor}\Label{main2}
Let $M\subset \C^N$ be a real-analytic hypersurface which does not
contain any nontrivial complex-analytic subvariety and $\S\subset
\C^{N'}$ the unit sphere, $N'> N\geq 2$. Then any CR map $f\colon
M\to \S$ of class $\6C^{N'-N+1}$ extends holomorphically to a neighborhood
of $M$ in $\C^N$.
\end{Cor}

Concerning the initial required regularity in Theorem \ref{main1},
we do not know whether the $\6C^{N'-N+1}$-smoothness assumption is
optimal. It is however known (see \cite{L85, G87, Ha, D}) that
$\6C^0$-smoothness (and even $\6C^{\alpha}$-smoothness,
$0<\alpha<1/6$; see \cite{S})
 is not enough to guarantee the holomorphic extendability of any map $f$ in Theorem \ref{main1}. On the other hand,
 when $M$
is a unit sphere too, it follows from the work of Huang \cite{H99}
 that Theorem \ref{main1} holds with only $\6C^2$-smoothness assumption on the map $f$
 provided $N'<2N-1$ (see also the
previous related work in \cite{CS1, Fa}). We should also mention
that for a map $f$ of class $\6C^{\infty}$, Theorem \ref{main1}
follows from the work of Pushnikov \cite{Pu1, Pu2} (whose proof,
however, is incomplete and contains a gap; see \cite{MMZ2}). In
the case where $M$ is assumed to be strongly pseudoconvex (at
every point) and $f$ is of class $\6C^{N'-1}$, Theorem \ref{main1}
was announced by Pinchuk in \cite{Pi} but a proof has not been
published.

Our strategy for the proof of Theorem \ref{main1} is similar to
that of the above mentioned result of Baouendi-Huang-Rothschild
\cite{BHR}. The first step is to prove that in the setting of
Theorem \ref{main1}, the map $f$ extends meromorphically to a
neighborhood of $M$ in $\C^N$. The desired holomorphic extension
of $f$ will then follow by applying a result due to Chiappari
\cite{Ch}. The main novelty of this paper consists in the proof of
the above mentioned meromorphic extension of $f$ regardless of the
codimension $N'-N$ (see Proposition \ref{car}). The ingredients of
the proof rely on another meromorphic extension result for a class
of CR ratios proved in \cite{MMZ1}, and an inductive dichotomy (in
Lemma \ref{work} below) showing that necessarily each component of
the map $f$ belongs to this class of ratios (see also \cite{Mir}
for a related argument in the one-codimensional case in another
context).

\section{Proof of Theorem \ref{main1}}
For basic concepts and notions about CR maps, we refer the reader
to \cite{Bog, BERbook}.
Theorem \ref{main1} is an immediate consequence of the following two results
and Lemma \ref{hor1} (i).

\begin{Pro}\Label{car}
Let $M\subset \C^N$ and $\S\subset \C^{N'}$ be as in Theorem {\rm
\ref{main1}} and $f\colon M\to \S$ a CR map of class
$\6C^{N'-N+1}$. Then $f$ extends meromorphically to a neighborhood
of $M$ in $\C^N$.
\end{Pro}

\begin{Thm}\Label{Chia}
Let $D\subset \C^N$ be a domain, $M\subset \partial D$ a
real-analytic hypersurface and $F\colon D\to \C^{N'}$ a
holomorphic map extending continuously up to $M$,  $N'\geq N\geq
2$. Assume that $F$ maps $D$ into the unit ball $\B\subset
\C^{N'}$, that $F|_{M}$ sends $M$ into the unit sphere $\S\subset
\C^{N'}$ and that $F|_M$ extends meromorphically to a neighborhood
of $M$ in $\C^N$. Then $F|_M$ extends holomorphically to a
neighborhood of $M$ in  $\C^N$.
\end{Thm}

The remainder of the paper is devoted to the proof
of Proposition \ref{car} while the second ingredient (Theorem \ref{Chia})
was already proved by Chiappari in \cite{Ch} (see
also \cite{CS2} for a previous version where the source manifold
is a sphere).
We should also mention that in the case $N'=N+1$ and $M$ is a real-algebraic hypersurface of
D'Angelo finite type, a result analogous to Proposition \ref{car} has been proved in
\cite{BHR}.

\subsection{Some preliminary facts}
We start by collecting a few  facts from the known literature.
Recall here that a real-analytic hypersurface $M\subset \C^N$ is
said to be {\em minimal} at a point $p_0\in M$ if there is no
complex-analytic hypersurface of $\C^N$ contained in $M$ through $p_0$ (see
e.g.\ \cite{Tr, Tu, BERbook}). The following lemma follows from
the same arguments as those of \cite[Lemma 6.2]{BHR}.

\begin{Lem}\Label{hor1}
Let $M\subset \C^N$ be a real-analytic hypersurface, $\S\subset
\C^{N'}$ the unit real sphere, $N'> N\geq 2$, and let $f\colon
M\to \S$ be a non-constant CR map of class $\6C^{2}$. Then,
\begin{enumerate}
\item[(i)] if $M$ is minimal at every point, it is
pseudoconvex and moreover, for every $p\in M$, $f$ extends near $p$ to a holomorphic map $F$ on
the pseudoconvex side of $M$ and $F$ maps this side into the unit ball $\B$ of $\C^{N'}$.
\item[(ii)] if $M$ does not contain any nontrivial
complex-analytic subvariety, it is pseudoconvex and the set of
strongly pseudoconvex points of $M$ is dense in $M$.
\end{enumerate}
\end{Lem}

Lemma \ref{hor1} (i) together with an elementary unique
continuation argument implies:

\begin{Lem}\Label{hor2}
In the situation of Proposition {\rm \ref{car}}, if $f$ is not
constant, the set of strongly pseudoconvex points of $M$ is dense
in $M$.
\end{Lem}

Lemma \ref{hor2} in conjunction with standard Hopf lemma type
arguments (see e.g.\  \cite{F, Hu94, BHR, EL}) implies in
particular:

\begin{Lem}\Label{huj}
In the situation of Proposition {\rm \ref{car}}, if $f$ is not
constant, it is an immersion at a generic point of $M$, i.e.\
there is a dense open subset $\Omega$ of $M$ such that for all
$p\in \Omega$ the differential $df(p)\colon \C T_pM\to \C
T_{f(p)}\S$ is injective. $($Here $\C T_pM$ $($resp.\ $\C
T_{f(p)}\S$ $)$ denotes the complexified tangent space of $M$
$($resp.\ of $\S$$)$ at $p$ $($resp.\ at $f(p)$$)$.$)$
\end{Lem}

Finally we state the hypersurface version of a meromorphic
extension result proved in \cite[Thereom 2.6]{MMZ1} that will be
very useful for the proof of Proposition \ref{car}.

\begin{Thm}\label{NME}
Let $W \subset\C^N$, $V\subset\C^k$ be open subsets, $M\subset W$
a connected real-analytic hypersurface, $G\colon M\to V$ a
continuous CR map and $\Phi,\Psi\colon {V}^*\times W \to\C$
holomorphic functions, where $V^*:=\{\1{\zeta}:\zeta \in V\}$.
Assume that $M$ is minimal at every point and that there exists a
nonempty open subset of $M$ where $\Psi(\1{G(z)},z)$ does not
vanish and where the quotient
\begin{equation}\Label{ratio}
H(z):=\frac{\Phi (\1{G(z)},z)}{\Psi (\1{G(z)},z)} \end{equation}
is CR. Then $\Psi (\1{G(z)},z)$ does not vanish on a dense open
subset $\widetilde{M}\subset M$ and $H$ extends from $\widetilde
M$ meromorphically to a neighborhood of $M$ in $\C^N$.
\end{Thm}

We should mention that a preliminary version of Theorem \ref{NME},
namely the case where $G$ is $\6C^{\infty}$ over $M$ and $H$ is CR
on a dense open subset of $M$, is contained in \cite{CPS99}. It
will be however important for the proof of Proposition \ref{car}
(and therefore of Theorem \ref{main1}) to have the continuous
version provided by Theorem \ref{NME}.

\subsection{Reflection identities and a linear system with
coefficients of meromorphic type} We start here the effective
proof of Proposition \ref{car}. It is enough to prove that $f$
extends meromorphically to a neighborhood of a given point $p_0\in
M$. Without loss of generality we may assume that $p_0=f(p_0)=0$,
that $f$ is not constant and moreover, that $\S$ is replaced
by the Heisenberg hypersurface i.e.\
\begin{equation}\Label{Heisenberg}
\F=\{(z',w')\in \C^{N'-1}\times \C:
\rho'(z',w',\1{z'},\1{w'}):={\tt Im}\, w'-|z'|^2=0\},
\end{equation}
with $|z'|^2=\sum_{j=1}^{N'-1}|z_j'|^2$,
$z'=(z_1',\ldots,z_{N'-1}')$. Shrinking $M$ if necessary near the
origin, we may assume that $M$ is connected and may also choose an
open connected neighborhood $U\subset \C^N$ of the origin such
that $M\subset U$ and such that we are given a family
$\LL=(\1L_1,\ldots,\1L_{N-1})$ of $(0,1)$ vector fields with
real-analytic coefficients in $U$ with
$\1L_1|_M,\ldots,\1L_{N-1}|_{M}$ spanning $T^{0,1}M$, the $(0,1)$
tangent bundle of $M$. In what follows, for $\gamma
=(\gamma_1,\ldots,\gamma_{N-1})\in \N^{N-1}$,
$\gamma_1+\ldots+\gamma_{N-1}=:|\gamma|\leq
N'-N+1$, $\1L^{\gamma}$ denotes the standard differential operator
$\1L_1^{\gamma_1}\ldots \1L_{N-1}^{\gamma_{N-1}}$ acting on $\6C^{N'-N+1}$-smooth
functions (or maps) on $M$. It will be
convenient to define the following classes of functions over $M$.

\begin{Def}\Label{EDP} For $M$, $f$, $\LL$ as above and $0\leq l\leq N'-N+1$,
we let $\6S_l$ denote the set of functions on $M$ that can be
written as a polynomial in $(\1L^{\gamma}\1f)_{|\gamma|\leq l}$
with coefficients that are real-analytic functions over $M$.
Denote also by $\6S_{-1}$ the set of all real-analytic functions
over $M$. For $-1\leq l\leq N'-N+1$, we also define $\6R_l$ to be
the set of all ratios of the form $a/b$ with $(a,b)\in (\6S_l)^2$
and $b\not \equiv 0$.
\end{Def}

Obviously the class $\6S_l$ forms a subring of those functions
over $M$ that are of class $\6C^{N'-N+1-l}$ and satisfies
$\6S_l\subset \6S_{l+1}$. Note also that since $M$ is everywhere
minimal, a function in $\6S_l$ can not vanish on any open subset
of $M$ unless it is identically zero (see e.g.\ \cite[Section
6]{MMZ1}). Therefore any ratio in $\6R_l$ is defined on a dense
open subset of $M$, and can not vanish on any open subset of this
dense set unless it is identically zero too. This fact will be
useful for the proof of Lemma \ref{work} below.

The main observation of this paper relies on the following lemma.

\begin{Lem}\Label{work}
In the situation of Proposition {\rm \ref{car}} and with the
above notation, shrinking $M$ around the origin if necessary, each
component of the map $f$ can be written as a ratio belonging to
$\6R_{N'-N+1}$.
\end{Lem}

\begin{proof} In the $(z',w')$ coordinates we split the map $f$
as follows $f=(\2f,g)=(\2f_1,\ldots,\2f_{N'-1},g)\in
\C^{N'-1}\times \C$. Since $f(M)\subset \F$, we have the relation
satisfied on $M$
\begin{equation}\Label{funda}
g-\1{g}=2i \left(\sum_{\nu=1}^{N'-1}|\2f_{\nu}|^2\right).
\end{equation}
By Lemma \ref{huj} the map $f$ is an immersion at a generic point
and in particular, for a generic point $p\in M$, $df(p)\colon
T^{0,1}_pM\to T_{f(p)}^{0,1}\F$ is injective. Pick such a point
$p_0\in M$. Then we may select $N-1$ components among those of
$\2f$, say $\2f_1,\ldots,\2f_{N-1}$, such that the matrix
$(\1L_j\1{\2f}_k)_{1\leq j,k\leq N-1}$ is invertible at $p_0$. The
above matrix $(\1L_j\1{\2f}_k)_{1\leq j,k\leq N-1}$ has
necessarily rank $N-1$ in a dense open subset of $M$ since its
determinant belongs to the class $\6S_1$ (see Definition
\ref{EDP}).  Now we split the map $\2f$ as follows
$\2f=(h,\psi)\in \C^{N-1}\times \C^{N'-N}$ where
$h=(h_1,\ldots,h_{N-1}):=(\2f_1,\ldots,\2f_{N-1})$ and $\psi
=(\psi_1,\ldots,\psi_{N'-N}):=(\2f_{N},\ldots,\2f_{N'-1})$, so
that $f=(h,\psi,g)$. Applying for each  $j=1,\ldots,N-1$ the
vector field $\1L_{j}$ to (\ref{funda}) and using the fact that
$f$ is CR, we obtain

\begin{equation}\Label{funda2}
-\1L_j\1g =2i \left( \sum_{\nu=1}^{N-1}h_{\nu}\cdot
\1L_j\1h_{\nu}+\sum_{\nu=1}^{N'-N}\psi_{\nu}\cdot
\1L_j\1\psi_{\nu}\right) \ {\rm on}\  M.
\end{equation}
Since the $\6C^{N'-N}$-smooth function given by $M\ni z\mapsto
{\rm det}\, \left(\1L_j\1{h}_{\nu}(z)\right)_{1\leq j,\nu\leq
N-1}$ does not vanish on a dense open subset $M_1$ of $M$, we may
apply Cramer's rule to (\ref{funda}) and (\ref{funda2}) to obtain
a relation of the form
\begin{equation}\Label{funda3}
\left(h,g\right)= r^{(1)}_{0}+\sum_{\nu=1}^{N'-N}\psi_{\nu}\cdot
r^{(1)}_{\nu}\quad {\rm on}\ M_1.
\end{equation}
Here for $\nu=0,\ldots,N'-N$, each $r^{(1)}_{\nu}$ is a $\C^N$-valued ratio
with components in the ring $\6R_1$, and
is defined on $M_1$. Since $f$ is of class at least $\6C^{2}$, we
may again apply for each $j=1,\ldots,N-1$ the vector field $\1L_j$ to
(\ref{funda3}) on $M_1$ to obtain :
\begin{equation}\Label{funda4}
0=\1L_jr^{(1)}_0+\sum_{\nu=1}^{N'-N}\psi_{\nu} \cdot \1L_j
r_{\nu}^{(1)}\quad {\rm on}\  M_1.
\end{equation}
Define the following $N(N-1)\times (N'-N)$ matrix with
$\6C^{N'-N-1}$-smooth coefficients over $M_1$ by setting
$$A^{(2)}:=\begin{pmatrix}
 \1L_1r^{(1)}_{1}&\ldots&\1L_{1}r_{N'-N}^{(1)}\\
\1L_2r^{(1)}_1&\ldots&\1L_{2}r_{N'-N}^{(1)}\\
\vdots&\ldots &\vdots\\
 \1L_{N-1}r^{(1)}_{1}&\ldots& \1L_{N-1}r^{(1)}_{N'-N}\end{pmatrix}, $$
 where each $r_{\nu}^{(1)}\in \C^N$ is viewed as a column
vector. Since any minor of $A^{(2)}$ has its coefficients in the
class $\6R_2$ the matrix $A^{(2)}$ achieves its maximal rank, that
we denote by $n^{(2)}$, on a dense open subset of $M_1$. Now
similarly to \cite{F, Hu94, BHR, Mir}, we come to a dichotomy
according to whether $n^{(2)}>0$ or $n^{(2)}=0$.

{\sc Case 1:} $n^{(2)}>0$. By taking a suitable minor of the
matrix $A^{(2)}$, which does not vanish on a dense open subset
$M_2$ of $M_1$, and using (\ref{funda4}) and Cramer's rule we see
that at least one of the components of $\psi$, say $\psi_{N'-N}$,
can be expressed as an affine combination of the following form:
\begin{equation}\Label{crook}
\psi_{N'-N}=c^{(2)}_0+\sum_{\nu=1}^{N'-N-1}\psi_{\nu}\cdot
c_{\nu}^{(2)}.
\end{equation}
Here for $\nu=0,\ldots,N'-N-1$, each $c^{(2)}_{\nu}$ is a
$\C^{N+1}$-valued ratio
with components in the ring $\6R_2$, and
is defined on $M_2$. By using
(\ref{funda3}) and (\ref{crook}), we may write
\begin{equation}\Label{funda5}
\left(h,\psi_{N'-N},g\right)=
r^{(2)}_{0}+\sum_{\nu=1}^{N'-N-1}\psi_{\nu}\cdot
r^{(2)}_{\nu},\quad {\rm on}\ M_2,
\end{equation}
where each $r^{(2)}_{\nu}$ is a $\C^{N+1}$-valued ratio with components in
$\6R_2$, and is defined on $M_2$.

{\sc Case 2:} $n^{(2)}=0$. We claim that also in this case, we may
reduce the starting system of equations (\ref{funda3}) on $M_1$ to
another system of the form (\ref{funda5}) that holds on a possibly
other dense subset of $M$, possibly after interchanging the
components of $\psi$ and shrinking $M$ near the origin. For this,
note that $n^{(2)}=0$ is equivalent to saying that each ratio
$r_{\nu}^{(1)}$, $\nu=1,\ldots,N'-N$, is a CR map on $M_1$, and
hence $r_0^{(1)}$ too in view of (\ref{funda4}). Since each
component of each ratio $r^{(1)}_{\nu}$ can clearly be written in
the form (\ref{ratio}) for a suitable $\Psi,\Phi$ and $G$
(see \cite[Lemma 6.1]{MMZ1} for more details), we
may apply Theorem \ref{NME} to conclude that each $r^{(1)}_{\nu}$
extends as a meromorphic mapping to a neighborhood of the origin
in $\C^N$. Let $u_{\nu}$ be a $\C^N$-valued holomorphic map and
$v_{\nu}$ a nonvanishing holomorphic function both defined near
the origin in $\C^N$ such that $u_{\nu}/v_{\nu}$ gives the
meromorphic extension of $r_{\nu}^{(1)}$ near $0$. Then after
shrinking $M$ near the origin if necessary, we may rewrite
(\ref{funda3}) as follows
\begin{equation}\Label{funda3bis}
(h(z),g(z))=\frac{u_0(z)}{v_0(z)}+\sum_{\nu=1}^{N'-N}\psi_{\nu}(z)\frac{u_{\nu}(z)}{v_{\nu}(z)},\
z\in \2M_1,
\end{equation}
on some suitable dense open subset $\2M_1$ of $M$. Following the
splitting used in (\ref{funda3bis}), for any $\nu=0,\ldots,N'-N$,
we write $u_{\nu}=(u_{\nu}',\hat u_{\nu})\in \C^{N-1}\times \C$.
Consider the family $(P_z)_{z\in \2M_1}$ of $N'-N$ dimensional
affine complex planes defined via the following parametrization
$\eta_z\colon \C^{N'-N}\to \C^{N'}$, $z\in \2M_1$,
\begin{equation}\Label{eta}
t=(t_1,\ldots,t_{N'-N})\mapsto \left(\frac{u'_0(z)}{v_0(z)}+
\sum_{\nu=1}^{N'-N}t_{\nu}\frac{u'_{\nu}(z)}{v_{\nu}(z)},
t_1,\ldots,t_{N'-N},\frac{\hat u_0(z)}{v_0(z)}+
\sum_{\nu=1}^{N'-N}t_{\nu}\frac{\hat
u_{\nu}(z)}{v_{\nu}(z)}\right). \end{equation}
 Since the Heisenberg hypersurface cannot contain any positive-dimensional affine complex
subspace, we have $P_z\not \subset \F$ for all $z\in \2M_1$ i.e.\
\begin{equation}\Label{aue}
\rho'(\eta_z(t),\1{\eta_z(t)})\not \equiv 0,\ {\rm for}\ z\in
\2M_1,\ t\in \C^{N'-N},
\end{equation}
where $\rho'$ is given by (\ref{Heisenberg}). It is clear that
$\rho'(\eta_z(t),\1{\eta_z(t)})$ may be written as follows
\begin{equation}\Label{omit}
\mu+\sum_{i=1}^{N'-N}\lambda_i\1{t_i}+\sum_{j=1}^{N'-N}
t_j\left(\sum_{k=1}^{N'-N}\xi_{jk}\1{t_k}+\sigma_{jk}\right),
\end{equation}
where $\mu$, $\lambda_i$, $\xi_{jk}$ and $\sigma_{jk}$ are ratios
in $\6R_{-1}$ and defined on $\2M_1$. (We have deliberately
omitted to write the dependence on $z$ for those ratios in
(\ref{omit}).) But since $f$ sends $M$ into $\F$, in view of
(\ref{eta}) and (\ref{funda3bis}), we have $\rho'(\eta_z(\psi
(z)),\1{\eta_z(\psi (z))})=0$ for $z\in \2M_1$ i.e.\
\begin{equation}\Label{omitbis}
\mu+\sum_{i=1}^{N'-N}\lambda_i\1{\psi_i}+\sum_{j=1}^{N'-N}\psi_j\left(\sum_{k=1}^{N'-N}(\xi_{jk}\1{\psi_k}+
\sigma_{jk})\right)=0\ {\rm on}\ \2M_1.
\end{equation}
Here again we come to a dichotomy.

{\sc Subcase 2.1:} there exists $j_0\in \{1,\ldots, N'-N\}$ such
that $\sum_{k=1}^{N'-N}(\xi_{j_0k}\1{\psi_k}+ \sigma_{j_0k})\not
\equiv 0$ (as a ratio in $\6R_0$). Then from (\ref{omitbis}) we
see that $\psi_{j_0}$ can be written as an affine combination,
with coefficients in $\6R_{0}\subset \6R_2$, of the other
components of $\psi$. Using (\ref{funda3}) we see that the claim
mentioned at the very beginning of {\sc Case 2} is proved.

{\sc Subcase 2.2:} for any $j\in \{1,\ldots, N'-N\}$,
\begin{equation}\Label{bus}
\sum_{k=1}^{N'-N}(\xi_{jk}\1{\psi_k}+ \sigma_{jk}) \equiv 0,\
{\rm on}\ \2M_1.
\end{equation}
If there exists $j_1,k_1\in \{1,\ldots,N'-N\}$ such that
$\xi_{k_1j_1}\not \equiv 0$, then after taking complex conjugates
of (\ref{bus}) , $\psi_{k_1}$ can be written as an affine
combination, with coefficients in $\6R_{-1}\subset \6R_2$, of the
other components of $\psi$. The desired claim therefore follows as
in Subcase 2.1. If not, then necessarily all $\xi_{jk}$ and
$\sum_{1\leq k\leq N'-N}\sigma_{jk}$ (by (\ref{bus})) are
identically zero. Equation (\ref{omitbis}) therefore leads to
$\mu+\sum_{1\leq i\leq N'-N}\lambda_i\1{\psi_i}=0$. Then there
exists necessarily $i_2\in \{1,,\ldots,N'-N\}$ such that
$\lambda_{i_2}\not \equiv 0$ since otherwise all quantities $\mu$,
$\lambda_i$, $\xi_{jk}$ and $\sum_{1\leq k\leq N'-N} \sigma_{jk}$
would be identically zero, which would contradict (\ref{aue}) in
view of (\ref{omit}). Hence here again we obtain that $\psi_{i_2}$
can be written as an affine combination, with coefficients in
$\6R_{-1}\subset \6R_2$, of the other components of $\psi$, which
completes, as explained before, the proof of the claim announced
at the beginning of {\sc Case 2}.

We are now able to finish the proof of Lemma \ref{work}. The
conclusion of the dichotomy studied in {\sc Cases} 1 and 2 states
that we can always  reduce a system of the form (\ref{funda3})
with coefficients in $\6R_1$ to another system of the form
(\ref{funda5}) with coefficients in $\6R_{2}$. If $N'=N+1$, we are
clearly done. If not, by using similar arguments as in the above
process, it is easy to show that the obtained system with
coefficients in $\6R_2$ can also be reduced to a system where
$N+2$ components of the map can be expressed as an affine
combination of the other ones with coefficients in $\6R_3$. (Note
that the only ingredients needed for such a process are Theorem
\ref{NME} and the fact that the Heisenberg hypersurface does not
contain any complex affine subspace of positive dimension.) By
pursuing this procedure $N'-N-2$ more times, we conclude that each
component of $f$ can be written as a ratio belonging to
$\6R_{N'-N+1}$. The proof of Lemma \ref{work} is complete.
\end{proof}

\subsection{Completion of the proof of Proposition {\rm \ref{car}}} By
Lemma \ref{work} we know that after shrinking $M$ near the
origin, each component of $f$ agrees on a dense open
subset of $M$ with a ratio in the class $\6R_{N'-N+1}$. Moreover since $M$ is minimal at $0$,
 $f$ extends to a holomorphic map $F$ defined on one side
 of $M$ near the origin \cite{BT84, Tr}. Therefore by using this extension $F$, it
is easy to see that any ratio in the class $\6R_{N'-N+1}$ can be written in
the form (\ref{ratio}) for
a suitable $\Psi$, $\Phi$ and $G$ (for further details
on that matter
 see \cite[Lemma 6.1]{MMZ1}). Since each
component of $f$
is CR, we conclude from Theorem \ref{NME} that each such component
extends meromorphically to a neighborhood of the origin in $\C^N$.
This completes the proof of Proposition \ref{car} and hence that of
Theorem \ref{main1}.

\begin{Rem}
In the proof of Theorem \ref{main1}, the assumption that $M$ does
not contain any complex-analytic hypersurface of $\C^N$ is
essential when applying several times Theorem \ref{NME}. We do not
know whether Theorem \ref{main1} holds without such an assumption.
(The main problem without such a hypothesis would be to prove the
meromorphic extension of the map $f$ near any non-minimal point.)
On the other hand, the assumption that $M$ be strongly
pseudoconvex at some point in Theorem \ref{main1} is used to know
that any non-constant CR map $f\colon M\to \S$ is a  local
embedding near a generic point of $M$. Without assuming the
existence of a strongly pseudoconvex point on $M$, the conclusion
of Theorem \ref{main1} still holds if $f$ is assumed to be of
class $\6C^{N'-r_f+1}$ where $r_f:={\rm Sup}\, \{{\rm Rank}\,
df(p)|_{T^{0,1}_pM}:p\in M\}$ (with a trivial modification in the
beginning of the proof of Lemma \ref{work}). This implies in
particular that Theorem \ref{main1} still holds for all CR maps of
class $\6C^{N'}$ without assuming the existence of any strongly
pseudoconvex point on $M$.
We leave the details to the reader.
\end{Rem}

We conclude this paper by giving examples of situations covered by
Theorem \ref{main1} and Corollary \ref{main2} where previously known results
do not apply.

\begin{Exa}
Given nonnegative integers $k,l$, $k\geq 1$, $l\geq 3$, and $h$ a
holomorphic function defined near the origin in $\C$ with
$h(0)=h'(0)=0$, define $M_k\subset \C^2$ to be the real-analytic
hypersurface through the origin given by
$$M_k:=\{(z_1,z_2)\in \C^2:{\tt Im}\, z_2=|z_1|^{2k}+|h(z_1)|^2\}.$$
Note that there are non-constant CR maps from $M_k$ into the
Heisenberg representation
\begin{equation}\Label{Heis}
\H:=\{(w_1,\ldots,w_l)\in \C^l: {\tt Im}\,
w_l=|w_1|^2+\ldots+|w_{l-1}|^2 \}
\end{equation}
 of the unit sphere in $\C^l$, namely the restriction to $M_k$ of the
holomorphic map
$$(z_1,z_2)\mapsto
(z_1^k,h(z_1),\underbrace{0,\ldots,0}_{\text{$l-3$ times}},z_2).$$
For $k=1$, $M_1$ is strongly pseudoconvex and the analyticity (at
every point of $M_1$) of $\6C^{l-1}$-smooth CR maps from $M_1$
into $\H$ follows from Corollary \ref{main2} but not from
\cite{F}. In the case $l=3$, the regularity of such maps does not
follow neither from \cite{BHR} if the function $h$ is chosen not
be algebraic.

For $k>1$, $M_k$ is weakly pseudoconvex of D'Angelo finite type
and the set of weakly pseudoconvex points is a real line. Here
again the analytic regularity (at every point of $M_k$) of
$\6C^{l-1}$-smooth CR maps from $M_k$ into $\H$ follows from
Corollary \ref{main2} but not from \cite{F}. Moreover, even if
$l=3$ and $h$ is chosen to be algebraic, the analyticity at the
origin of such maps does not follow from \cite{BHR} since the type
of $M_k$ at $0$ is clearly greater or equal to 4 (and the required
initial order of smoothness for the maps in \cite{BHR} equals this
type).
\end{Exa}

\begin{Exa}\Label{add}
Let $M\subset \C^3$ be the real-algebraic hypersurface given by
$M:=\{(z_1,z_2,z_3)\in \C^3:{\tt Im}\, z_3=|z_1z_2|^2+|z_2|^2\}$.
Then $M$ obviously contains a complex line through the origin and
the set of strongly pseudoconvex points of $M$ is dense in $M$.
(Note also that $M$ is not essentially finite at the origin in the
sense of \cite{BJT}.) For any nonnegative integer $l\geq 4$, the
restriction to $M$ of the following holomorphic map
$$(z_1,z_2,z_3)\mapsto
(z_1z_2,z_2,\underbrace{0,\ldots,0}_{\text{$l-3$ times}},z_3)\in
\C^l$$ maps $M$ into $\H$ (given by (\ref{Heis})), and is near a
generic point a local embedding. The analyticity (at every point
of $M$) of $\6C^{l-2}$-smooth CR maps from $M$ into $\H$ follows
from Theorem \ref{main1} but neither from \cite{F} nor from
Corollary \ref{main2}.
\end{Exa}

\section*{Acknowledgements}
This paper was written while I was visiting UCSD. I would like to
thank this institution for its hospitality and M.S. Baouendi, M.
Derridj, P. Ebenfelt, B. Lamel, F. Meylan, L.P. Rothschild and D.
Zaitsev for useful discussions and comments on this paper. I am
also grateful to the referee for several helpful remarks.

\end{document}